\documentclass[11pt,bezier]{article}
\usepackage{amsmath,graphicx,amssymb,amsfonts,color}

\font\bg=cmbx10 scaled\magstep1 
\font\small=cmr8

\newtheorem{newlemma}{{\bf Lemma}}

\newenvironment{lema}{\begin{newlemma}{\hspace{-0.5
em}{\bf.}}}{\end{newlemma}}

\newtheorem{newteorem}{{\bf Theorem}}

\newenvironment{teorem}{\begin{newteorem}{\hspace{-0.5
em}{\bf.}}}{\end{newteorem}}

\newtheorem{newkorolari}{{\bf Corollary}}

\newenvironment{korolari}{\begin{newkorolari}{\hspace{-0.5
em}{\bf.}}}{\end{newkorolari}}

\newtheorem{newexample}{{\bf Example}}

\newenvironment{example}{\begin{newexample}{\hspace{-0.5
em}{\bf.}}}{\end{newexample}}
\newtheorem{newkonjek}{{\bf Conjecture}}

\newtheorem{newdefine}{{\bf Definition}}

\begin{document}
\tolerance=10000 \baselineskip18truept
\newbox\thebox
\global\setbox\thebox=\vbox to 0.2truecm{\hsize
0.15truecm\noindent\hfill}
\def\boxit#1{\vbox{\hrule\hbox{\vrule\kern0pt
     \vbox{\kern0pt#1\kern0pt}\kern0pt\vrule}\hrule}}
\def\qed{\lower0.1cm\hbox{\noindent \boxit{\copy\thebox}}\bigskip}
\def\ss{\smallskip}
\def\ms{\medskip}
\def\nt{\noindent}

\centerline{\Large \bf Independent sets of some graphs associated}  \centerline{\Large \bf
to commutative rings}
 \vspace{.3cm}


\bigskip
\baselineskip12truept \centerline{\bf Saeid Alikhani$^a$ and Saeed Mirvakili$^b$
}\baselineskip20truept \centerline{$^a$Department of Mathematics, Yazd
University, 89195-741, Yazd,
Iran} \centerline{ $^{b}$Department of Mathematics, Payame Noor University, 19395-4697 Tehran, I.R. Iran}

 \nt\rule{16cm}{0.1mm}

\nt{\bg ABSTRACT}

\medskip

\nt Let $G=(V,E)$ be a simple graph. A set $S\subseteq V$ is
independent set of $G$,  if no two vertices of $S$ are adjacent.
The  independence number $\alpha(G)$ is the size of a maximum
independent set in the graph. 
Let $R$ be a commutative ring with 
nonzero identity and  $I$ an ideal of
$R$. The zero-divisor graph of $R$, denoted by $\Gamma(R)$, is
an undirected graph whose vertices are the nonzero zero-divisors
of $R$ and two distinct vertices $x$ and $y$ are adjacent if and
only if $xy = 0$. Also the ideal-based zero-divisor graph of $R$, denoted by $\Gamma_I(R)$, is the graph which
vertices are the set $\{x\in R\backslash I | xy\in I  \quad\mbox{for some} \quad  y\in R\backslash I\}$	 and two distinct vertices $x$
and $y$ are adjacent if and only if $xy \in I$.
In this paper we study the independent sets and the independence number of
 $\Gamma(R)$ and $\Gamma_I(R)$.

\nt {\bf Keywords:}  {\small Independent set; Independence
number; Zero-divisor graph, Ideal.}

\nt{\bf 2010 Mathematics Subject Classification:} {\small 05C69, 13A99.}

\nt\rule{16cm}{0.1mm}

\section{Introduction}

\nt A simple graph $G=(V,E)$ is a finite nonempty set $V(G)$ of
objects called vertices together with a (possibly empty) set
$E(G)$ of unordered pairs of distinct vertices of $G$ called
edges. The concept of zero-divisor graph of a commutative ring
with identity was introduced by Beck in \cite{beck}  and has been
studied in \cite{akbari1,akbari2,anderson1,anderson2,axtel}. Redmond in \cite{redmond} has extended
this concept to any arbitrary ring.  Let $R$ be a commutative ring
with $1$. The zero-divisor graph of $R$, denoted $\Gamma(R)$, is
an undirected graph whose vertices are the nonzero zero-divisors
of $R$ and two distinct vertices $x$ and $y$ are adjacent if and
only if $xy = 0$. Thus $\Gamma(R)$ is an empty graph if and only
if $R$ is an integral domain.

\nt The concept of dominating set in zero-divisor graph has implicitly
been studied in \cite{naderi} and \cite{mojdeh}. Throughout this article,
all rings are commutative with identity $1 \neq 0$. For a subset
$A$ of a ring $R$, we let $A^* = A\setminus \{0\}$. We will denote
the rings of integers modulo $n$, the integers, and a finite field
with $q$ elements by $\mathbb{Z}_n, \mathbb{Z}$ and $F_q$, respectively. For
commutative ring theory, see \cite{ati,kaplan}.

\nt An  independent set of a graph $G$ is a set of vertices where
no two vertices are adjacent. The  independence number $\alpha(G)$
is the size of a maximum independent set in the graph. An
independent set with cardinality
 $\alpha(G)$ is called a $\alpha$-set (\cite{alikhani,gutman,hoede}). 

\nt A graph $G$ is called bipartite if its vertex set can be
partitioned into $X$ and $Y$ such that every edge of $G$ has one
endpoint in $X$ and other endpoint in $Y$. 
A complete
$r$-partite graph is one whose vertex set can be partitioned into
$r$ subsets so that no edge has both ends in any one subset and
each vertex of a partite set is joined to every vertex of the
another partite sets. We denote a complete bipartite graph by $K_{m,n}$.
The graph $K_{1,n}$ is called a star
graph, and a bistar graph is a graph generated by two $K_{1,n}$ graphs, where their centers
are joined.
For a nontrivial connected graph $G$ and a
pair vertices $u$ and $v$ of $G$, the distance $d(u,v)$ between
$u$ and $v$ is the length of a shortest path from $u$ to $v$ in
$G$. The girth of a graph $G$, containing a cycle, is the smallest
size of the length of the cycles of $G$ and is denoted by $gr(G)$.
If $G$ has no cycles, we define the girth of $G$ to be infinite. A
graph in which each pair of distinct vertices is joined by an edge
is called a complete graph $K_n$ on $n$ vertices.  For a graph $G$,
a complete subgraph of $G$ is called a clique. The clique number,
$\omega(G)$, is the greatest integer $n \geq 1$ such that
$K_n\subseteq G$, and $\omega(G)$ is infinite if $K_n \subseteq G$
for all $n \geq  1$, see \cite{west}.

\nt Dominating sets and domination number of zero-divisor
graphs and ideal-based zero-divisor graphs investigated  in \cite{mojdeh}. Similar to \cite{mojdeh}, in this paper,  we study
independent sets and independence numbers  of zero-divisor
graphs and ideal-based zero-divisor graphs. In Section 2 we review some preliminary results related
to independence number of a graph. In Section 3, we study the independence number of zero-divisor graphs associated to commutative
rings. In Section 4, investigate the independence number of an ideal based zero-divisor graph of $R$. Finally in Section 5, we list tables for graphs associated to small commutative ring $R$ with 1, and write independence, domination and clique number of
$\Gamma(R)$.

\section{Preliminary results}

There are several classes of graphs whose independent sets and
independence  numbers are clear. We state some of them in the
following Lemma, which  their proofs are straightforward.

\begin{lema}{\rm (\cite{west})}\label{west}
\begin{enumerate}
\item[(i)] $\alpha(K_n)=1$.

\item[(ii)] Let $G$ be a complete $r$-partite graph $(r\geq 2)$
with partite sets $V_1, ...,V_r$. If $|V_i|\geq 2$ for $1\leq i
\leq r$, then $\alpha(G) = max_{1\leq i \leq r}|V_i|$.

\item[(iii)] $\alpha(K_{1,n}) = n$ for a star graph $K_{1,n}$.

\item[(iv)] The independence  number of a bistar graph is $2n$.

\item[(v)]  Let $C_n, P_n$ be a cycle and a path with $n$
vertices, respectively. Then $\alpha(P_n)=\lfloor
\frac{n+1}{2}\rfloor$ and $\alpha(C_n)=\lfloor
\frac{n}{2}\rfloor$.

\end{enumerate}
\end{lema}

\begin{korolari}
Let $F_1$ and $F_2$ be finite fields with $|F_1^*|= m $ and
$|F_2^*|= n$. Then

\begin{enumerate}
\item[(i)]  $\alpha(\Gamma(F_1\times F_2))=\max\{m,n\}.$

\item[(ii)] $\alpha(\Gamma(F_1 \times \mathbb{Z}_4))=\max\{2m,3\}.$

\end{enumerate}

\end{korolari}

\nt{\bf Proof.} \begin{enumerate}
\item[(i)] The graph $\Gamma(F_1\times F_2)$ is bipartite (\cite{anderson1}) and we have the result by Lemma \ref{west} (ii).

\item[(ii)] We have  $Z^*(F_1 \times \mathbb{Z}_4)=\{(x,y)|x\in\ F_1^*, y=0,2\}\cup\{(0,y)|y=1,2,3\}$.

 If $F_1=\Bbb Z_2$ then  $\{(0,y)|y=1,2,3\}$ is  a maximum independent set in the graph and so $\alpha(\Gamma(F_1 \times\Bbb  Z_4))=3$.
If $F_1\neq\Bbb Z_2$ then $\{(x,y)|x\in\ F_1^*, y=0,2\}$ is  a maximum independent set in the graph and so $\alpha(\Gamma(F_1 \times\Bbb  Z_4))=2m$. Therefore $\alpha(\Gamma(F_1 \times \mathbb{Z}_4))=\max\{2m,3\}.$\quad\qed
\end{enumerate}

\section{Independence number of a zero-divisor graph}

\nt We start  this section with the following lemma:

\begin{lema} Let $R$ be a ring and $r\geq 3$. If $\Gamma(R)$ is a  $r$-partite graph with parts
$V_1,\ldots,V_r$, then $\alpha(\Gamma(R))=max|V_i|$.
\end{lema}

\nt Note that,  for any prime number $p$ and any positive integer
$n$, there exists a finite ring $R$ whose zero-divisor graph
$\Gamma(R)$ is a complete $p^n$-partite graph. For example, if
$\Gamma(R)$ is a finite field with $p^n$ elements, then $R
=F_{p^n}[x,y]/(xy,y^2-x)$ is the desired ring.

\nt {\bf Remark.} It is easy to see that a graph $G$ has
independence  number equal to $1$, if for every $x,y\in Z(R)^*$, $xy=0$, this means $\Gamma(R)$ is a complete graph.


\ms
\nt We need the following theorem:

\begin{teorem} {\rm(\cite{anderson2})}
 \label{int} If $R$ is a commutative ring which is not an integral domain,
then exactly one of the following holds:

\begin{enumerate}
\item[(i)] $\Gamma(R)$ has a cycle of length $3$ or $4$ (i.e.,
$gr(R) \leq 4)$;

\item[(ii)]  $\Gamma(R)$ is a star graph; or

\item[(iii)] $\Gamma(R)$ is the zero-divisor graph of $R\cong
\mathbb{Z}_2 \times \mathbb{Z}_4$ or $R\cong \mathbb{Z}_2 \times
\mathbb{Z}[X]/(X^2)$.

\end{enumerate}
\end{teorem}

\nt By Theorem \ref{int} we have the following theorem:

\begin{teorem}
If $\Gamma(R)$ has no cycles, then $\alpha(\Gamma(R))$ is either
$|Z^*(R)|-1$ or $3$.
\end{teorem}

\begin{teorem}

\begin{enumerate}
\item[(i)] Let $R$ be a finite ring. If $\Gamma(R)$ is a regular
graph of degree $r$, then $\alpha(\Gamma(R))$ is either 1 or $r$.

\item[(ii)]  Let $R$ be a finite decomposable ring. If $\Gamma(R)$
is a Hamiltonian graph, then $\alpha(\Gamma(R))=\frac{|Z^*(R)|}{2}$.

\item[(iii)] Let $R$ be a finite principal ideal ring and not decomposable. If
$\Gamma(R)$ is Hamiltonian, then  $\alpha(\Gamma(R))=1$
\end{enumerate}
\end{teorem}
\nt{\bf Proof.}
\begin{enumerate}
\item[(i)]
  Since $\Gamma(R)$ is a regular graph of degree
$r$, $\Gamma(R)$ is a complete graph $K_{r+1}$ or a complete bipartite
graph $K_{r,r}$. Consequently, $\alpha(\Gamma(R))$ is either 1 or
$r$.

\item[(ii)] In this case $\Gamma(R)$ is $K_{n,n}$ for some natural
number $n$. So, $\alpha(\Gamma(R))=n$.

\item[(iii)] If $R$ is not decomposable then in this case $\Gamma(R)$ is  a complete graph. Therefore we have the
result.\quad\qed

\end{enumerate}
\begin{korolari} The graph $\Gamma(\mathbb{Z}_n)$ is a Hamiltonian
graph if and only if $\alpha(\Gamma(\mathbb{Z}_n))=1$.
\end{korolari}
\nt{\bf Proof.} By Corollary 2 of \cite{akbari2}, we know that the graph
$\Gamma(\mathbb{Z}_n)$ is a Hamiltonian graph if and only if
$n=p^2$, where $p$ is a prime larger than $3$ and
$\Gamma(\mathbb{Z}_n)$ is isomorphic to $K_{p-1}$. Thus, we have
the result.\quad\qed

\nt Here we state a notation which is useful for the study of $\alpha$-sets and the independence number of graphs associated to commutative rings.

\nt Let $R=F_1\times\ldots\times F_n$, where  $F_i$ is an integral domain, for every $i$, and $|F_i|\geq |F_{i+1}|.$ We consider $E_{i_1\ldots i_k}$ as the following set: $$E_{i_1\ldots i_k}=\{(x_1,\ldots,x_n)\in R| \forall i\in \{i_1,\ldots,i_k\}, x_i\neq 0\ {\rm and\ } \forall i\not\in \{i_1,\ldots,i_k\}, x_i=0\}.$$
By this notation we have $|E_{i_1\ldots i_k}|=|F_{i_1}^*||F_{i_2}^*|\ldots|F_{i_k}^*|.$ In the rest of this paper we use the following equation:
\[
n_{i_1}n_{i_2}...n_{i_k}=|E_{i_1i_2...i_k}|.
\]

\begin{teorem} Suppose that for a fixed integer $n\geq 2$, $R = F_1 \times \cdots \times F_n$, where $F_i$ is an integral
domain for each $i = 1,\ldots,n$. We have
\begin{enumerate}
\item[(i)] $\alpha(\Gamma(R))= \infty $ if   one of $F_i$ is infinity,

\item[(ii)] \[\alpha(\Gamma(R))\geq  \left(\sum_{2\leq i_2\leq \ldots \leq i_{\lfloor\frac{k-1}{2}\rfloor}\leq n}n_1n_{i_2}\ldots n_{i_{\lfloor\frac{k-1}{2}\rfloor}}\right)+\sum_{\lfloor\frac{k-1}{2}\rfloor+1}^{n-1}\left(\sum_{1\leq i_1\leq \ldots \leq i_l\leq n} n_{i_1}\ldots n_{i_l}\right).
 \]

\end{enumerate}
\end{teorem}
\nt{\bf Proof.} $(i)$ We can suppose that  $|F_1|$ is infinity. So $S=\{(x,0,\ldots,0)|x\in R_1^*\}$ is an independent set and therefore $\alpha(\Gamma(R))=\infty.$

$(ii)$ Let $|F_1|\geq |F_2|\geq\ldots\geq|F_n|.$  It is easy to see that
\[A= \left(\displaystyle\bigcup_{2\leq i_2\leq \ldots \leq i_{\lfloor\frac{k-1}{2}\rfloor}\leq n}E_{1i_2\ldots i_{\lfloor\frac{k-1}{2}\rfloor}}\right) \bigcup \left(  \bigcup_{\lfloor\frac{k-1}{2}\rfloor+1}^{n-1}\left(\bigcup_{1\leq i_1\leq \ldots \leq i_l\leq n} E_{i_1\ldots i_l}\right)\right)
\]
is an independent set of $\Gamma(R)$.   So
 \[\alpha(\Gamma(R))\geq |A|=\sum_{\lfloor\frac{k-1}{2}\rfloor+1}^{n-1}\left(\sum_{1\leq i_1\leq \ldots \leq i_l\leq n} n_{i_1}\ldots n_{i_l}\right)+ \left(\sum_{2\leq i_2\leq \ldots \leq i_{\lfloor\frac{k-1}{2}\rfloor}\leq n}n_1n_{i_2}\ldots n_{i_{\lfloor\frac{k-1}{2}\rfloor}}\right)\quad\qed
 \]

\begin{teorem}
Suppose that  $n_1\geq n_2\geq n_3$ and $|F_i^*|=n_i$ for  $i=1,2,3$. If $R=F_1\times F_2\times F_3$, then $$\alpha(\Gamma(R))= n_1n_2+n_1n_3+\max\{n_1,n_2n_3\}.$$
\end{teorem}\label{f3}
\nt{\bf Proof.} It is not difficult to see that one of the following sets is a $\alpha$-set in the zero-divisor graph of $R$:
\begin{enumerate}
\item[] $A_1=E_{12}\cup E_{13} \cup E_{23}$,
\item[] $A_2=E_{12}\cup E_{13}\cup E_{1}.$
\end{enumerate}
So $\alpha(\Gamma(R))= \max\{|A_1|,|A_2|\}=n_1n_2+n_1n_3+\max\{n_1,n_2n_3\}.$ \quad\qed

\nt Let us to state two examples for the above theorem:
\begin{example}
Let $R=\Bbb Z_5\times \Bbb Z_2\times\Bbb Z_2$. Here $A_2=E_{12}\cup E_{13}\cup E_{1}$ is a $\alpha$-set  of graph $\Gamma(R)$ and so $\alpha(\Gamma(R))=n_1n_2+n_1n_3+n_1=9.$
\end{example}
\begin{example}
Let $R=\Bbb Z_7\times \Bbb Z_5\times\Bbb Z_5$. Here  $A_1=E_{12}\cup E_{13} \cup E_{23}$ is  a $\alpha$-set  and $\alpha(\Gamma(R))=n_1n_2+n_1n_3+n_2n_3=64.$
\end{example}
\begin{teorem}\label{f4}
Suppose that  $n_1\geq n_2\geq n_3\geq n_4$ and $|F_i^*|=n_i$ for  $i=1,2,3,4$. Let $R=F_1\times F_2\times F_3\times F_4$.
 \begin{enumerate}
\item[(i)]  If $n_1\geq n_2n_3n_4$, then $\alpha(\Gamma(R))=n_1(n_2n_3+n_2n_4+n_3n_4+n_2+n_3+n_4+1).$

\item[(ii)] If $n_1\leq n_2n_3n_4$ and $n_1n_4\geq n_2n_3$, then $\alpha(\Gamma(R))=n_1(n_2n_3+n_2n_4+n_3n_4+n_2+n_3+n_4)+n_2n_3n_4.$

\item[(iii)] If  $n_1n_4\leq n_2n_3$, then $\alpha(\Gamma(R))=n_1(n_2n_3+n_2n_4+n_3n_4+n_2+n_3)+n_2n_3+n_2n_3n_4.$
    \end{enumerate}
\end{teorem}
\nt{\bf Proof.} Since $n_1\geq n_2\geq n_3\geq n_4$, it is easy to check that one of the following sets is a $\alpha$-set of  the graph
$\Gamma(R)$:
\begin{enumerate}
\item[] $I_1=  E_{123}\cup E_{124}\cup E_{134}\cup E_{12}\cup E_{13}\cup E_{14}\cup E_{1},$
\item[] $I_2=  E_{123}\cup E_{124}\cup E_{134}\cup E_{12}\cup E_{13}\cup E_{14}\cup E_{234},$
\item[] $I_3=E_{123}\cup E_{124}\cup E_{134}\cup E_{12}\cup E_{13}\cup E_{23}\cup E_{234},$
\end{enumerate}
\begin{enumerate}
\item[(i)] Suppose that  $n_1\geq n_2n_3n_4$. Obviously we have $n_1n_4\geq n_2n_3$. So in this case   $I_1$  is a $\alpha$-set in the graph. Therefore $\alpha(\Gamma(R))=n_1(n_2n_3+n_2n_4+n_3n_4+n_2+n_3+n_4+1).$
\item[(ii)] If $n_1\leq n_2n_3n_4$ and $n_1n_4\geq n_2n_3$, then in this case $I_2$ is  a $\alpha$-set in the graph. Therefore $\alpha(\Gamma(R))=n_1(n_2n_3+n_2n_4+n_3n_4+n_2+n_3+n_4)+n_2n_3n_4.$

\item[(iii)] If  $n_1n_4\leq n_2n_3$ then $n_1\leq n_2n_3n_4$ and so in this case $I_3$ is  a $\alpha$-set in the graph. So $\alpha(\Gamma(R))=n_1(n_2n_3+n_2n_4+n_3n_4+n_2+n_3)+n_2n_3+n_2n_3n_4.$\quad\qed
    \end{enumerate}
    \nt The following corollary is an immediate consequence of Theorem \ref{f4}.

\begin{korolari}
   Suppose that  $n_1\geq n_2\geq n_3\geq n_4$ and $|F_i^*|=n_i$ for  $i=1,2,3,4$.
   If  $R=F_1\times F_2\times F_3\times F_4$, then \[\alpha(\Gamma(R))=n_1(n_2n_3+n_2n_4+n_3n_4+n_2+n_3)+\max\{n_1+n_1n_4,n_2n_3+n_2n_3n_4,n_1n_4+n_2n_3n_4\}.\]
\end{korolari}

\nt Here we state   some examples for Theorem \ref{f4}.
\begin{example}
Let $R=\Bbb Z_5\times \Bbb Z_2\times\Bbb Z_2\times\Bbb Z_2$. The set $I_1$  in Theorem \ref{f4}, is a $\alpha$-set in the graph and so $\alpha(\Gamma(R))=28.$
\end{example}
\begin{example}
Let $R=\Bbb Z_5\times \Bbb Z_3\times\Bbb Z_3\times\Bbb Z_3$.  The set $I_2$  in Theorem \ref{f3}, is a $\alpha$-set in the graph and so $\alpha(\Gamma(R))=80.$
\end{example}
\begin{example}
Let $R=\Bbb Z_5\times \Bbb Z_5\times\Bbb Z_3\times\Bbb Z_2$. The set $I_3$ in Theorem \ref{f3}, is a $\alpha$-set in the graph and so $\alpha(\Gamma(R))=88.$
\end{example}

\begin{teorem}\label{main}
Suppose that  $|F_i^*|=n_i$, where $n_i\geq n_j$ and $i\geq j$ for $i,j=1,\ldots,5$. Let $R=F_1\times \ldots\times F_5$, and
 $t=n_1(\displaystyle\sum_{2\leq i<j<k\leq 5}n_in_jn_k)+n_1(\displaystyle\sum_{\stackrel{2\leq i<j\leq 5}{(i,j)\neq(4,5)}}n_in_j)$. We have
\begin{enumerate}
\item[(i)]  If $n_1\geq n_2n_3n_4n_5$, then $\alpha(\Gamma(R))=t+n_1(n_4n_5+n_2+n_3+n_4+n_5+1).$

\item[(ii)] If $n_2n_3\geq n_1n_4n_5$,  then $\alpha(\Gamma(R))=t+n_2(n_3n_4n_5+n_3n_4+n_3n_5+n_1+n_3)+n_1n_3.$

\item[(iii)] If  $ n_1n_5\geq n_2n_3n_4$, then $\alpha(\Gamma(R))=t+n_1(n_4n_5+n_2+n_3+n_4+n_5)+n_2n_3n_4n_5.$

\item[(iv)]  If $ n_1n_5\leq n_2n_3n_4$ and $ n_1n_4\geq n_2n_3n_5$, then $\alpha(\Gamma(R))=t+n_1(n_4n_5+n_2+n_3+n_4)+n_2(n_3n_4n_5+n_3n_4).$

\item[(v)]  If $ n_1n_4\leq n_2n_3n_5$ and $ n_1n_3\geq n_2n_4n_5$, then $\alpha(\Gamma(R))=t+n_1(n_4n_5+n_2+n_3)+n_2(n_3n_4n_5+n_3n_4+n_3n_5).$

\item[(vi)]  If $ n_1n_3\leq n_2n_4n_5$ and $ n_1n_2\geq n_3n_4n_5$, then $\alpha(\Gamma(R))=t+n_1(n_4n_5+n_2)+n_2(n_3n_4n_5+n_3n_4+n_3n_5+n_4n_5).$

\item[(vii)]  If $ n_1n_2\leq n_3n_4n_5$, then $\alpha(\Gamma(R))=t+(n_1+n_3)n_4n_5+n_2(n_3n_4n_5+n_3n_4+n_3n_5+n_4n_5).$

    \end{enumerate}
\end{teorem}
\nt{\bf proof.} We put $A=(\displaystyle\bigcup_{2\leq i<j<k\leq 5}E_{1ijk})\bigcup(\displaystyle\bigcup_{\stackrel{2\leq i<j\leq 5}{(i,j)\neq(4,5)}}E_{1ij}).$ Consider the sets $A_i$ and $B_i$ for $i=1,\ldots,6$ as shown in the following table.
\[
\begin{tabular}{|l|l|l|l|l|l|l|}
  \hline
         &  $  i= 1 $     &   $i=  2$       &   $ i= 3 $    &$ i=   4 $     &   $i=   5$  &$i= 6$ \\
\hline$ A_i$ &$ E_1  $    & $E_{23}$  & $E_{12}$  & $E_{13}$  & $E_{14}$&$E_{15}$ \\
\hline $B_i $& $E_{2345} $& $E_{145}$ & $E_{345}$ & $E_{245}$ & $E_{235}$ & $E_{234}$\\
  \hline
\end{tabular}
\]
We have:
\begin{enumerate}
\item[(i)]  If $n_1\geq n_2n_3n_4n_5$, then by the above table $|A_1|\geq|B_1|$ and this implies $|B_2|\geq |A_2|$ and for $i=3,4,5,6$, $|A_i|\geq |B_i|$. So $A\cup A_1\cup B_2\cup A_3\cup A_4\cup A_5\cup A_6$ has the size of a maximum independent set in the graph and $\alpha(\Gamma(R))=t+n_1(n_4n_5+n_2+n_3+n_4+n_5+1).$

\item[(ii)] If $n_2n_3\geq n_1n_4n_5$  then $|A_2|\geq|B_2|$ and this implies $|B_1|\geq |A_1|$,$|A_3|\geq|B_3|$,$|A_4|\geq|B_4|$,$|B_5|\geq|A_5|$ and $|B_6|\geq|A_6|$, so $A\cup B_1\cup A_2\cup A_3\cup A_4\cup B_5\cup B_6$ has the size of a maximum independent set in the graph and $\alpha(\Gamma(R))=t+n_2(n_3n_4n_5+n_3n_4+n_3n_5+n_1+n_3)+n_1n_3.$

\item[(iii)] If  $ n_1n_5\geq n_2n_3n_4$ and $n_1\leq n_2n_3n_4n_5$ then $|A_6|\geq|B_6|$ and $|B_1|\geq|A_1|,$ now $|B_2|\geq |A_2|$ and for $i=3,4,5$, $|A_i|\geq |B_i|$, so $A\cup B_1\cup B_2\cup A_3\cup A_4\cup A_5\cup A_6$ has the size of a maximum independent set in the graph and  $\alpha(\Gamma(R))=t+n_1(n_4n_5+n_2+n_3+n_4+n_5)+n_2n_3n_4n_5.$

\item[(iv)]  If $ n_1n_5\leq n_2n_3n_4$ and $ n_1n_4\geq n_2n_3n_5$ then $|B_6|\geq|A_6|$ and $|A_5|\geq|B_5|,$ now $|B_1|\geq |A_1|$, $|B_2|\geq |A_2|$ and for $i=3,4$, $|A_i|\geq |B_i|$, so $A\cup B_1\cup B_2\cup A_3\cup A_4\cup A_5\cup B_6$ has the size of a maximum independent set in the graph and   $\alpha(\Gamma(R))=t+n_1(n_4n_5+n_2+n_3+n_4)+n_2(n_3n_4n_5+n_3n_4).$

\item[(v)]  If $ n_1n_4\leq n_2n_3n_5$ and $ n_1n_3\geq n_2n_4n_5$ then $|B_5|\geq|A_5|$ and $|A_4|\geq|B_4|,$ therefore  $|A_3|\geq |B_3|$ and for $i=1,2,6$, $|B_i|\geq |A_i|$, so $A\cup B_1\cup B_2\cup A_3\cup A_4\cup B_5\cup B_6$ has the size of a maximum independent set in the graph and    $\alpha(\Gamma(R))=t+n_1(n_4n_5+n_2+n_3)+n_2(n_3n_4n_5+n_3n_4+n_3n_5).$

\item[(vi)]  If $ n_1n_3\leq n_2n_4n_5$ and $ n_1n_2\geq n_3n_4n_5$ then $|B_4|\geq|A_4|$ and $|A_3|\geq|B_3|,$ so   for $i=1,2,5,6$, $|B_i|\geq |A_i|$, hence $A\cup B_1\cup B_2\cup A_3\cup B_4\cup B_5\cup B_6$ has the size of a maximum independent set in the graph and     $\alpha(\Gamma(R))=t+n_1(n_4n_5+n_2)+n_2(n_3n_4n_5+n_3n_4+n_3n_5+n_4n_5).$

\item[(vii)]  If $ n_1n_2\leq n_3n_4n_5$ then $|B_3|\geq|A_3|$ and   for $i=1,2,4,5,6$, $|B_i|\geq |A_i|$, hence $A\cup B_1\cup B_2\cup B_3\cup B_4\cup B_5\cup B_6$ has the size of a maximum independent set in the graph and     $\alpha(\Gamma(R))=t+(n_1+n_3)n_4n_5+n_2(n_3n_4n_5+n_3n_4+n_3n_5+n_4n_5).$\quad\qed
    \end{enumerate}

\begin{korolari}\label{f5}
Let $R=F_1\times \ldots\times F_5$, $|F_i^*|=n_i$  and $n_i\geq n_j$, where  $i,j=1,\ldots,5$ and $i\geq j.$ Then $$\alpha(\Gamma(R))=n_1(\sum_{2\leq i<j<k\leq 5}n_in_jn_k)+n_1(\sum_{\stackrel{2\leq i<j\leq 5}{(i,j)\neq(4,5)}}n_in_j)+\max_i\{\Delta_i\},$$
where
\[\begin{array}{l}
\Delta_1=n_1(n_4n_5+n_2+n_3+n_4+n_5+1)\\ \Delta_2=n_2(n_3n_4n_5+n_3n_4+n_3n_5+n_1+n_3)+n_1n_3\\
\Delta_3=n_1(n_4n_5+n_2+n_3+n_4+n_5)+n_2n_3n_4n_5\\ \Delta_4=n_1(n_4n_5+n_2+n_3+n_4)+n_2(n_3n_4n_5+n_3n_4) \\
\Delta_5=n_1(n_4n_5+n_2+n_3)+n_2(n_3n_4n_5+n_3n_4+n_3n_5)\\ \Delta_6=n_1(n_4n_5+n_2)+n_2(n_3n_4n_5+n_3n_4+n_3n_5+n_4n_5)\\
\Delta_7=(n_1+n_3)n_4n_5+n_2(n_3n_4n_5+n_3n_4+n_3n_5+n_4n_5)
\end{array}\]

\end{korolari}

\begin{example}
\begin{enumerate}
\item[(i)] Let $R=\Bbb Z_5\times \Bbb Z_2\times\Bbb Z_2\times\Bbb Z_2\times\Bbb Z_2$. Then in Theorem \ref{main}, $\alpha(\Gamma(R))=t+\Delta_1,$
\item[(ii)] Let $R=\Bbb Z_5\times \Bbb Z_5\times\Bbb Z_5\times\Bbb Z_2\times\Bbb Z_2$. Then in Theorem \ref{main}, $\alpha(\Gamma(R))=t+\Delta_2,$
\item[(iii)] Let $R=\Bbb Z_7\times \Bbb Z_3\times\Bbb Z_3\times\Bbb Z_3\times\Bbb Z_3$. Then in Theorem \ref{main}, $\alpha(\Gamma(R))=t+\Delta_3,$
\item[(iv)] Let $R=\Bbb Z_7\times \Bbb Z_3\times\Bbb Z_3\times\Bbb Z_3\times\Bbb Z_2$. Then in Theorem \ref{main}, $\alpha(\Gamma(R))=t+\Delta_4,$
\item[(v)] Let $R=\Bbb Z_5\times \Bbb Z_5\times\Bbb Z_5\times\Bbb Z_2\times\Bbb Z_2$. Then in Theorem \ref{main}, $\alpha(\Gamma(R))=t+\Delta_5,$
\item[(vi)] Let $R=\Bbb Z_7\times \Bbb Z_7\times\Bbb Z_3\times\Bbb Z_2\times\Bbb Z_2$. Then in Theorem \ref{main}, $\alpha(\Gamma(R))=t+\Delta_6,$
\item[(vii)] Let $R=\Bbb Z_3\times \Bbb Z_3\times\Bbb Z_3\times\Bbb Z_3\times\Bbb Z_3$. Then in Theorem \ref{main}, $\alpha(\Gamma(R))=t+\Delta_7.$
\end{enumerate}
\end{example}
\begin{teorem}
Let $(R,\texttt{m})$ be a finite local ring and $\texttt{m}\neq\{0\}$.
\begin{enumerate}
 \item[(i)] If $\texttt{m}^2=\{0\}$, then $\alpha(\Gamma(R))=1.$
 \item[(ii)] If $\texttt{m}^2\neq\{0\}$, then $2\leq\alpha(\Gamma(R))\leq |Z^*(R)|-|Ann(Z(R))^*|.$
 \end{enumerate}
\end{teorem}
\nt{\bf Proof.} If R is a finite local ring, then the Jacobson radical of R
equals $Z(R)$ and $Z(R)=\texttt{m}$. Thus $Z(R)$ is a nilpotent ideal and since $R$ is not a field, then
$Ann(Z(R)) \neq \{0\}.$ Moreover, each element of $Ann(Z(R))$
 is adjacent to each other vertex
of $Z^*(R)$.

$(i)$ If $\texttt{m}^2=\{0\}$ then $Ann(Z(R))=Z^*(R)$ and $\Gamma(R)$ is a complete graph.

$(ii)$ If $\texttt{m}^2\neq\{0\}$, then every element of $Ann(Z(R))^*$
 is adjacent to each other vertex of $Z^*(R)$ and this implies $2\leq\alpha(\Gamma(R))\leq |Z^*(R)|-|Ann(Z(R))^*|.$\quad\qed

\begin{example}
Let $R=\Bbb Z_{p^3}$ then $Z^*(R)=\{pk|(p,k)=1\}\cup\{p^2k|(p^2,k)=1\}.$ We have $Ann(Z(R))^*=\{p^2k|(p^2,k)=1\}$ and $\{pk|(p,k)=1\}$ is an independent set in the
$\Gamma(R)$ of maximum size. So $\alpha(\Gamma(R))=|\{pk|(p,k)=1\}|=|Z^*(R)|-|Ann(Z(R))^*|.$
\end{example}


\section{The independence number of an ideal-based zero-divisor graph}

\nt Suppose that  $R$ is  a commutative ring with nonzero identity, and  $I$ is an ideal of
$R$. The ideal-based zero-divisor graph of $R$, denoted by $\Gamma_I(R)$, is the graph which
vertices are the set $\{x\in R\backslash I | xy\in I  \quad\mbox{for some}\quad  y\in R\backslash I\}$	 and two distinct vertices $x$
and $y$ are adjacent if and only if $xy \in I,$ see \cite{Redmond3}. In the case $I = 0$, $\Gamma_0(R)$ is denoted
by $\Gamma(R)$. Also, $\Gamma_I(R)$ is empty if and only if $I$ is prime. Note that Proposition 2.2(b)
of \cite{Redmond3} is equivalent to saying $\Gamma_I(R)=\emptyset$  if and only if $R/I$ is an integral domain.
That is, $\Gamma_I(R)=\emptyset$ if and only if $\Gamma(R/I)=\emptyset$.

\nt In this section, we study  the relationship between   independence numbers
of $\Gamma_I(R)$ and $\Gamma(R/I)$.

\nt This naturally raises
the question: If $R$ is  a commutative ring with ideal $I$,  whether $\alpha(\Gamma_I(R))$ is equal to $\alpha(\Gamma(R/I))?$ We show that the answer is negative in general.

\begin{lema} Let $m$ be a composite natural number and $p$  a prime number. Then
\[
\alpha(\Gamma_{m\Bbb Z}(\Bbb Z))=\left\{
\begin{array}{lr}
{\displaystyle \alpha(\Bbb Z/m\Bbb Z)=1 };&
\quad\mbox{if $m=p^2$,}\\[15pt]
{\displaystyle \infty};&
\quad\mbox{otherwise.}
\end{array}
\right.
\]
Note that for the second case $\alpha(\Gamma_{m\Bbb Z}(\Bbb Z))=\infty$ and $\alpha(\Bbb Z/m\Bbb Z)<\infty.$

\end{lema}
\nt{\bf Proof.}  If $m=p^2$ then for every $x\in \Gamma_{m\Bbb Z}(\Bbb Z)$ we have $x=pk$, where $(p,k)=1.$ So  $x,y\in \Gamma_{m\Bbb Z}(\Bbb Z)$ are adjacent in $\Gamma_I(R)$ and $\Gamma_I(R)$ is a complete graph. Also $\Bbb Z/m\Bbb Z\cong\Bbb Z_{p^2}$ and $\Gamma(\Bbb Z/m\Bbb Z)$ is a complete graph. Therefore in this case  $\alpha(\Gamma_{m\Bbb Z}(\Bbb Z))=\alpha(\Bbb Z/m\Bbb Z)=1$.

\nt  Now suppose that there isn't any prime number $p$ such that $m= p^2$. Then we have $m=p^in,$ where $p$ is prime number, $n\neq1$ and $(n,p)=1$, or $m=p^l,$ $p$ is prime and $l\geq 3$.

\nt If $m=p^l$ then $S=\{kp|(k,p)=1\}$ is an independent set and therefore $\alpha(\Gamma_{m\Bbb Z}(\Bbb Z))=\infty.$

\nt If $m=p^in$ then $S=\{kp|(k,p)=1 \ {\rm and\ } n|k\}$ is an independent set and therefore $\alpha(\Gamma_{m\Bbb Z}(\Bbb Z))=\infty.$
But, we have $\Bbb Z/m\Bbb Z$ is a finite ring and $\alpha(\Gamma(\Bbb Z/m\Bbb Z))$ is finite.\quad\qed

\nt Now we state the following results of \cite{Redmond3}.

\begin{lema} \label{lemma16}{\rm(\cite{Redmond3})}
 Let $I$ be an ideal of a ring $R$, and $x, y$ be in $R\backslash I$.
Then:\begin{enumerate}
\item[(i)]  If $x + I$ is adjacent to $y + I$ in $\Gamma(R/I)$, then $x$ is adjacent to $y$ in $\Gamma_I(R)$;
\item[(ii)]  If $x$ is adjacent to $y$ in $\Gamma_I(R)$ and $x + I \neq y + I$, then $x + I$ is adjacent to $y + I$ in
$\Gamma(R/I)$;
\item[(iii)] If $x$ is adjacent to y in $\Gamma_I(R)$ and $x + I = y + I$, then $x^2, y^2 \in I$.
\end{enumerate}
\end{lema}

\begin{lema} {\rm(\cite{Redmond3})}
 If $x$ and $y$ are (distinct) adjacent vertices in $\Gamma_I(R)$,
then all (distinct) elements $x + I$ and $y + I$ are adjacent in $\Gamma_I(R)$. If $x^2 \in I$, then all the
distinct elements of $x + I$ are adjacent in $\Gamma_I(R)$.
\end{lema}

\begin{teorem} Let $S$ be a nonempty subset of $R\backslash I$. If $S + I = \{s + I|s\in S	\}$ is an independent set of $\Gamma(R/I)$,
then $S$ is an independent set  of $\Gamma_I(R)$.
\end{teorem}
\nt {\bf Proof}. Let $S$ be a  nonempty subset of $R\backslash I$ and $S + I = \{s + I|s\in S	\}$  an independent set of $\Gamma(R/I)$. If $x,y\in S $,
then $x+I$ and $y+I$ are not adjacent in $\Gamma(R/I)$ and by Lemma \ref{lemma16}(i),  $x$ and $y$ are not adjacent  in $\Gamma_I(R)$.\quad\qed

\nt The following corollary is an immediate consequence of the above theorem:

\begin{korolari}
 $\alpha(\Gamma(R/I))\leq \alpha(\Gamma_I(R)).$
\end{korolari}

\begin{teorem}\label{AT}
Let $S+I$ be an independent set  with cardinality $\alpha(\Gamma(R/I))$ and $A=\{s+I\in S+I| s^2+ I=I\}$. Then $\alpha(\Gamma_I(R))=|A|+|I|(\alpha(\Gamma(R/I))-|A|).$
\end{teorem}

\nt{\bf Proof.} Suppose that $s\in S$, $x\in s+I$ and $y\in s+I$. If $s^2\in I$ then   $x\in s+I$ and $y\in s+I$ are adjacent vertices in $\Gamma_I(R).$ If $s^2\not\in I$ then $x\in s+I$ and $y\in s+I$ are not adjacent in $\Gamma_I(R).$ Therefore $T=\{s|s^2\in I\}\cup \{s+i|i\in I, s^2\not\in I\}$ is an independent set  with maximum cardinality.\quad\qed

\begin{korolari}
 $\alpha(\Gamma(R/I))\leq \alpha(\Gamma_I(R))\leq |I|\alpha(\Gamma(R/I))$
\end{korolari}

\begin{korolari}
If $S$ is an independent set  with cardinality $\alpha(\Gamma_I(R))$, and $s^2\in I$ for every $s\in S$, then $\alpha(\Gamma_I(R))=\alpha(\Gamma(R/I)).$
\end{korolari}

\begin{korolari}
If $S$ is an independent set  with cardinality $\alpha(\Gamma_I(R))$, and $s^2\not\in I$ for every $s\in S$, then $\alpha(\Gamma_I(R))=|I|\alpha(\Gamma(R/I)).$
\end{korolari}
\nt We state the following examples for above corollaries:

\begin{example}
Let $R=\Bbb Z_6\times \Bbb Z_3$ and $I=0\times \Bbb Z_3$ be an ideal of $R.$ Then it easy to see that $\Gamma_I(R)=\{(2,0),(2,1),(2,2),(3,0),(3,1),(3,2),(4,0),(4,1),(4,2)\}$ and $\Gamma(R/I)=\{(2,0)+I,(3,0)+I,(4,0)+I\}.$  The set $T=\{(2,0),(2,1),(2,2),(4,0),(4,1),(4,2)\}$ is an independent set of $\Gamma_I(R)$ and so
$\alpha(\Gamma_I(R))=6$. On the other hand $S+I=\{(2,0)+I,(4,0)+I\}$ is an independent set of $\Gamma(R/I)$ and
$\alpha(\Gamma(R/I))=2.$ Therefore $\alpha(\Gamma_I(R))=|I|\alpha(\Gamma(R/I)).$
\end{example}
\begin{example}
Let $R=\Bbb Z_{16}$ and $I=4\Bbb Z_{16}$. Then  $\Gamma_I(R)=\{2,6,10,14\}$ and $\Gamma(R/I)=\{2+I\}.$ Then $T=\{2\}$ is an independent set of $\Gamma_I(R)$ and
$\alpha(\Gamma_I(R))=1$. On the other hand $S+I=\{2+I\}$ is an independent set of $\Gamma(R/I)$ and
$\alpha(\Gamma(R/I))=1.$ So we have $\alpha(\Gamma_I(R))=\alpha(\Gamma(R/I)).$
\end{example}

\begin{example}
Let $R=\Bbb Z_{16}\times \Bbb Z_3$ and $I=0\times \Bbb Z_3$ be an ideal of $R.$ Then it easy to see that $\Gamma_I(R)=\{(x,y)|x=2,4,\ldots,14,y=0,1,2\}$ and $\Gamma(R/I)=\{(x,0)+I|x=2,4,\ldots,14\}.$ Then $T=\{(x,y)|x=2,6,10,14,y=0,1,2\}\cup\{(4,0)\}$ is an independent set of $\Gamma_I(R)$ and so
$\alpha(\Gamma_I(R))=13$. On the other hand $S+I=\{(x,0)+I|x=2,4,6,10,14\}$ is an independent set of $\Gamma(R/I)$ and
$\alpha(\Gamma(R/I))=5.$ Let $A$ be the set defined  in Theorem \ref{AT}, then $A=\{4\}$. So we have   $\alpha(\Gamma_I(R))=13=1+3(5-1)=|A|+|I|(\alpha(\Gamma(R/I))-|A|).$
\end{example}

\section{ Independence, domination and clique number of graphs associated to small finite commutative rings}
\nt In this section  similar to \cite{redmond2}, we list the tables for graphs associated to commutative ring $R$, and write independence, domination and clique number of $\Gamma(R)$.
Note that the tables for $n =|V\Gamma| = 1, 2, 3,  4$ can be found in \cite{anderson1}. The results for $n = 5$ can be found in \cite{Redmond3}.
In \cite{redmond2}, all graphs on $6, 7, \ldots, 14$ vertices which can be realized as the zero-divisor graphs of a commutative rings
with 1, and the list of all rings (up to isomorphism) which produce these graphs, has given.

\begin{figure}[h]
\hspace{4.8cm}
\includegraphics[width=4cm, height=2cm]{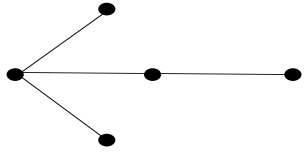}
\caption{\label{figure1} Graph for $\mathbb{Z}_2\times \mathbb{Z}_4$ and $\mathbb{Z}_2\times\mathbb{Z}_2[X]/(X^2)$ }
\end{figure}

\begin{figure}[h]
\hspace{4.8cm}
\includegraphics[width=4cm, height=2cm]{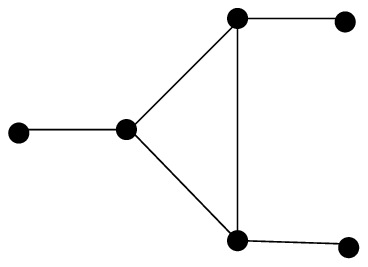}
\caption{\label{figure2} Graph for $\mathbb{Z}_2\times \mathbb{Z}_2\times\mathbb{Z}_2$ }
\end{figure}

\begin{figure}[h]
\hspace{4cm}
\includegraphics[width=5cm, height=2cm]{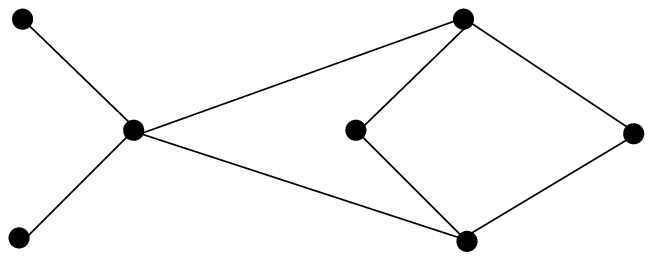}
\caption{\label{figure3} Graph for $\mathbb{Z}_3\times \mathbb{Z}_4$ and $\mathbb{Z}_3 \times\mathbb{Z}_2[X]/(X^2)$ }
\end{figure}

\begin{figure}[h]
\hspace{4.9cm}
\includegraphics[width=4.5cm, height=2.3cm]{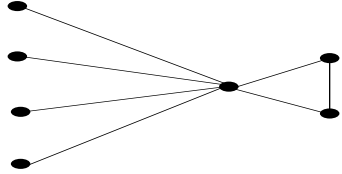}
\hspace{4cm}
\caption{\label{figure4} Graph for $\mathbb{Z}_4, \mathbb{Z}_2[X]/(X^4), \mathbb{Z}_4[X]/(X^2+2), \mathbb{Z}_4[X]/(X^2+3X)$ and
$\mathbb{Z}_4[X]/(X^3-2,2X^2,2X)$ }
\end{figure}

\begin{figure}[h]
\hspace{5cm}
\includegraphics[width=5cm, height=3cm]{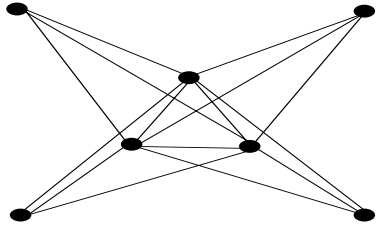}
\caption{\label{figure5} Graph for $\mathbb{Z}_2[X,Y]/(X^3,XY,Y^2), \mathbb{Z}_8[X]/(2X,X^2)$ and  $\mathbb{Z}_4[X]/(X^3,2X^2,2X)$ }
\end{figure}

\begin{figure}[h]
\hspace{5cm}
\includegraphics[width=4cm, height=2.2cm]{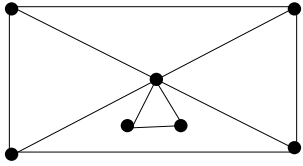}
\caption{\label{figure6} Graph for $\mathbb{Z}_4[X]/(X^2+2X), \mathbb{Z}_8[X]/(2X,X^2+4), \mathbb{Z}_2[X,Y]/(X^2,Y^2-XY)$
and $\mathbb{Z}_4[X]/(X^2,Y^2-XY,XY-2,2X,2Y)$ }
\end{figure}

\[
\begin{tabular}{|c|c|c|c|c|c|c|}
  \hline
  Vertices & R &$ |R|$ & Graph &   $\alpha(\Gamma(R))$ & $\gamma(\Gamma(R))$&$\omega(\Gamma(R)) $\\
 \hline 1 & $\Bbb Z_4 $& 4 & $K_1$ & 1 & 1&1 \\
 \hline   & $\Bbb Z_2[X]/(X^2) $ & 4& $K_1$ & 1 & 1 &1\\
  \hline
\end{tabular}
\]
\[
\begin{tabular}{|c|c|c|c|c|c|c|}
  \hline
  Vertices & R &$ |R|$ & Graph &   $\alpha(\Gamma(R))$ & $\gamma(\Gamma(R))$&$\omega(\Gamma(R)) $\\
 \hline 2 & $\Bbb Z_9 $& 9 & $K_2$ & 1 & 1&2 \\
 \hline   & $\Bbb Z_2\times\Bbb Z_2 $ & 9 & $K_2$ & 1 & 1&2\\
 \hline   & $\Bbb Z_3[X]/(X^2) $& 9 & $K_2$ & 1 & 1&2\\
  \hline
\end{tabular}
\]

\[
\begin{tabular}{|c|c|c|c|c|c|c|}
  \hline
  Vertices & R &$ |R|$ & Graph &   $\alpha(\Gamma(R))$ & $\gamma(\Gamma(R))$&$\omega(\Gamma(R)) $\\
 \hline 3 & $\Bbb Z_6 $& 6 & $K_{1,2}$ & 2 & 1&2 \\
 \hline   & $\Bbb Z_8$ & 8& $K_{1,2}$ & 2 & 1 &2\\
 \hline   & $\Bbb Z_2[X]/(X^3) $& 8 & $K_{1,2}$ & 2 & 1&2\\
  \hline  & $\Bbb Z_4[X]/(2X,X^2-2)$ & 8 &$ K_{1,2}$ & 2 & 1&2 \\
  \hline  & $\Bbb Z_2[X,Y]/(X,Y)^2$ & 8 & $K_3$ & 1 &1 &3\\
  \hline  & $\Bbb Z_4[X]/(2,X)^2$ & 8 & $K_3$ & 1 & 1&3 \\
 \hline   & $\Bbb F_4[X]/(X^2)$ & 16 & $K_3$ & 1 & 1 &3\\
 \hline   & $\Bbb Z_4[X]/(X^2+X+1)$ & 16 & $K_3$ & 1 & 1 &3\\
  \hline
\end{tabular}
\]

\begin{figure}[h]
\hspace{5.7cm}
\includegraphics[width=2.3cm, height=1.8cm]{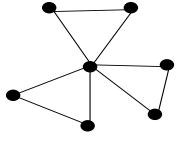}
\caption{\label{figure7} Graph for $\mathbb{Z}_4[X,Y]/(X^2,Y^2,XY-2,2X,2Y),\mathbb{Z}_2[X,Y]/(X^2,Y^2)$ and $\mathbb{Z}_4[X]/(X^2)$ }
\end{figure}

\begin{figure}[h]
\hspace{4.7cm}
\includegraphics[width=5cm, height=1.7cm]{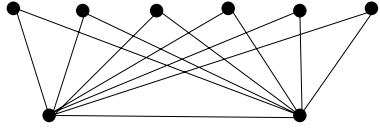}
\caption{\label{figure8} Graph for $\mathbb{Z}_4[X]/(X^3-X^2-2,2X^2,2X)$ }
\end{figure}

\begin{figure}[h]
\hspace{4.7cm}
\includegraphics[width=4.5cm, height=2cm]{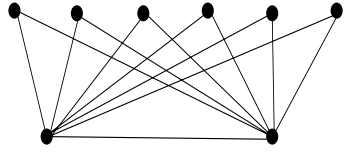}
\caption{\label{figure9} Graph for $\mathbb{Z}_9[X]/(3X,X^2-3), \mathbb{Z}_9[X]/(3X,X^2-6)$ and $\mathbb{Z}_3[X]/(X^3) $ }
\end{figure}

\begin{figure}[h]
\hspace{4.7cm}
\includegraphics[width=4.5cm, height=2cm]{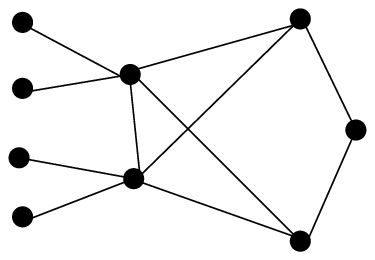}
\caption{\label{figure10} Graph for $\mathbb{Z}_2\times \mathbb{Z}_2\times \mathbb{Z}_2 $ }
\end{figure}

\begin{figure}[h]
\hspace{4.7cm}
\includegraphics[width=5cm, height=2cm]{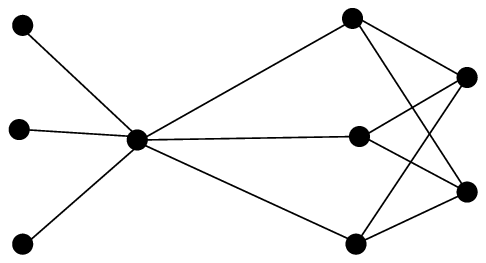}
\caption{\label{figure11} Graph for $\mathbb{Z}_4\times\Bbb F_4, \mathbb{Z}_2[X]/(X^2)\times F_4$ }
\end{figure}

\begin{figure}[h]
\hspace{4.7cm}
\includegraphics[width=5cm, height=2.9cm]{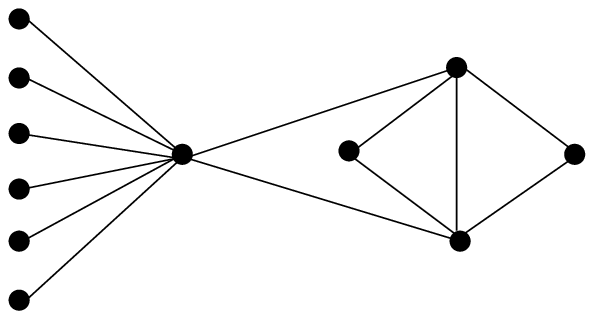}
\caption{\label{figure12} Graph for $\mathbb{Z}_2\times \mathbb{Z}_9$ and $\mathbb{Z}_2\times \mathbb{Z}_3[X]/(X^2) $ }
\end{figure}

\begin{figure}[h]
\hspace{4.7cm}
\includegraphics[width=4cm, height=2.2cm]{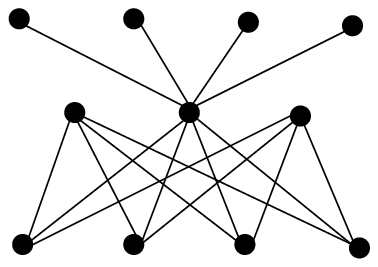}
\caption{\label{figure13} Graph for $\mathbb{Z}_5\times \mathbb{Z}_4$ and $\mathbb{Z}_5\times \mathbb{Z}_2[X]/(X^2) $ }
\end{figure}

\[
\begin{tabular}{|c|c|c|c|c|c|c|}
  \hline
  Vertices & R &$ |R|$ & Graph &   $\alpha(\Gamma(R))$ & $\gamma(\Gamma(R))$&$\omega(\Gamma(R)) $\\
 \hline 4 & $\Bbb Z_2\times\Bbb  F_4 $& 8 & $K_{1,3}$ & 3 & 1&2 \\
 \hline   & $\Bbb Z_3\times\Bbb Z_3 $ & 9& $K_{2,2}$ & 2 & 2 &2\\
 \hline   & $\Bbb Z_25 $& 25 & $K_4$ & 1& 1&4 \\
  \hline  & $\Bbb Z_5[X]/(X^2)$ & 25 &$ K_4$ & 1 & 1&4 \\
    \hline
\end{tabular}
\]
\[
\begin{tabular}{|c|c|c|c|c|c|c|}
  \hline
  Vertices & R &$ |R|$ & Graph &   $\alpha(\Gamma(R))$ & $\gamma(\Gamma(R))$&$\omega(\Gamma(R)) $\\
 \hline 5 & $\Bbb Z_2\times\Bbb  Z_5 $& 10 & $K_{1,4}$ & 4 & 1&2 \\
 \hline   & $\Bbb Z_3\times\Bbb F_4 $ & 12& $K_{2,3}$ & 3 & 2 &2\\
 \hline   & $\Bbb Z_2\times\Bbb Z_4 $& 8 & Fig. 1 & 3& 2&2 \\
  \hline  & $\Bbb Z_2\times Z_2[X]/(X^2)$ & 8 &  Fig. 1 & 2 & 1&2 \\
    \hline
\end{tabular}
\]

\begin{figure}[h]
\hspace{4.7cm}
\includegraphics[width=4.5cm, height=2.2cm]{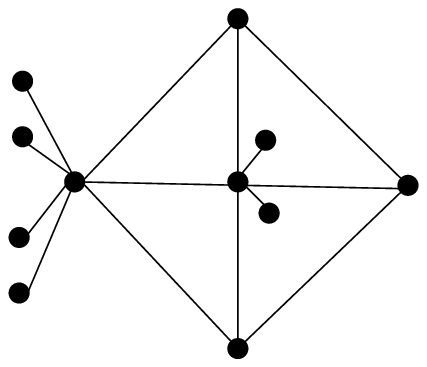}
\caption{\label{figure14} Graph for $\mathbb{Z}_2\times \mathbb{Z}_8$ and $\mathbb{Z}_2\times \mathbb{Z}_2[X]/(X^3)$ and
$\mathbb{Z}_2\times \mathbb{Z}_4[X]/(2X,X^2-2) $ }
\end{figure}

\begin{figure}[h]
\hspace{4.7cm}
\includegraphics[width=5cm, height=2.9cm]{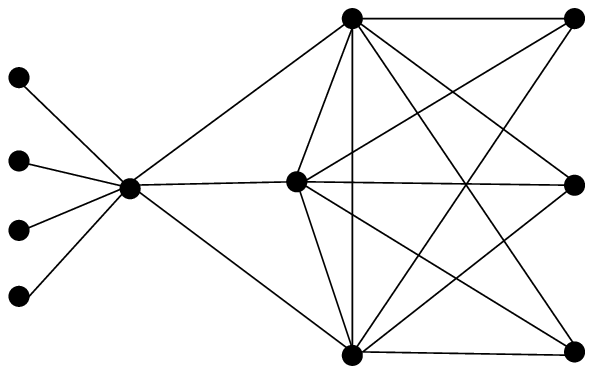}
\caption{\label{figure15} Graph for $\mathbb{Z}_2\times \mathbb{Z}_2[X,Y]/(X,Y)^2$ and $\mathbb{Z}_2\times \mathbb{Z}_4[X]/(2,X)^2$ }
\end{figure}

\begin{figure}[h]
\hspace{4.7cm}
\includegraphics[width=5cm, height=2.9cm]{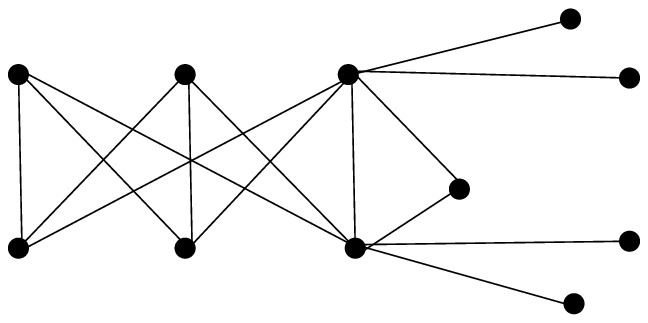}
\caption{\label{figure16} Graph for $\mathbb{Z}_4\times \mathbb{Z}_4, \mathbb{Z}_4\times \mathbb{Z}_2[X]/(X^2)$ and $\mathbb{Z}_2[X]/(X^2)\times \mathbb{Z}_2[X]/(X^2)$ }
\end{figure}

\[
\begin{tabular}{|c|c|c|c|c|c|c|}
  \hline
  Vertices & R &$ |R|$ & Graph &   $\alpha(\Gamma(R))$ & $\gamma(\Gamma(R))$&$\omega(\Gamma(R)) $\\
 \hline 6 & $\Bbb Z_3\times\Bbb  Z_5 $& 15 & $K_{2,4}$ & 4 & 2&2 \\
 \hline   & $\Bbb F_4\times\Bbb F_4 $ & 16& $K_{3,3}$ & 3 & 2 &2\\
 \hline   & $\Bbb Z_2\times\Bbb Z_2\times\Bbb Z_2 $& 8 & Fig. 2 & 3& 3&3 \\
  \hline  & $\Bbb Z_{49}$ & 49 &  $K_6$ & 1 & 1&6 \\
   \hline  & $\Bbb  Z_7[X]/(X^2)$ & 49 & $K_6$ & 1 & 1&6 \\
    \hline
\end{tabular}
\]


\begin{figure}[h]
\begin{minipage}{7.5cm}
\hspace{1.7cm}
\includegraphics[width=4.5cm,height=3.5cm]{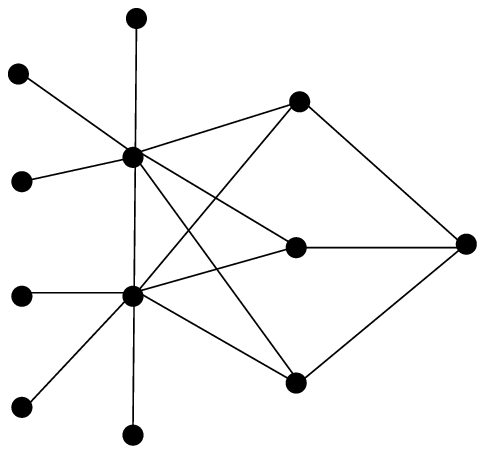}
\caption{\label{figure17} Graph for $\mathbb{Z}_2\times \mathbb{Z}_2 \times F_4 $}
\end{minipage}
\begin{minipage}{7.5cm}
\hspace{1.7cm}
\includegraphics[width=4.5cm,height=3.5cm]{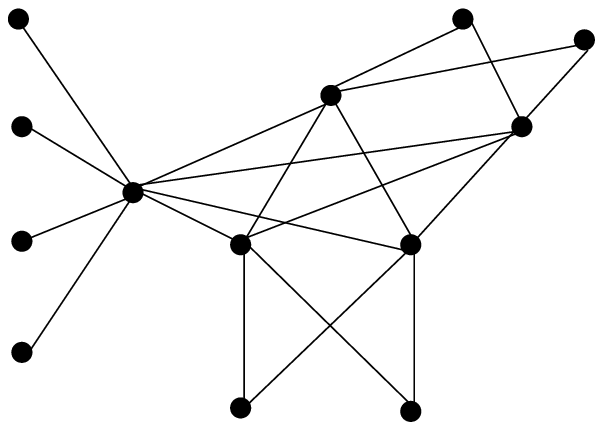}
\caption{\label{figure18} Graph for $\mathbb{Z}_2\times \mathbb{Z}_3 \times \mathbb{Z}_3 $.}
\end{minipage}
\end{figure}


\[
\begin{tabular}{|c|c|c|c|c|c|c|}
  \hline
  Vertices & R &$ |R|$ & Graph &   $\alpha(\Gamma(R))$ & $\gamma(\Gamma(R))$&$\omega(\Gamma(R)) $\\
 \hline 7 & $\Bbb Z_2\times\Bbb  Z_7 $& 14 & $K_{1,6}$ & 6 & 1&2 \\
 \hline   & $\Bbb F_4\times\Bbb Z_5 $ & 10& $K_{3,4}$ & 4 & 2 &2\\
 \hline   & $\Bbb Z_3\times\Bbb Z_4 $& 12 & Fig. 3 & 4& 2&2 \\
  \hline  & $\Bbb Z_3\times\Bbb Z_2[X]/(X^2)$ & 12 &  Fig. 3& 4& 2&2 \\
   \hline  & $\Bbb  Z_{16}$ & 16 & Fig. 4 & 5 & 1&3 \\
   \hline   & $\Bbb Z_2[X]/(X^4) $ & 16& Fig. 4 & 5 & 1&3\\
 \hline   & $\Bbb Z_4[X]/(X^2+2) $& 16 & Fig. 4 & 5 & 1&3 \\
  \hline  & $\Bbb Z_4[X]/(X^2+3X)$ & 16 &  Fig. 4 & 5 & 1&3 \\
   \hline  & $\Bbb  Z_4[X]/(X^3-2,2X^2,2X)$ & 16 & Fig. 4 &5 & 1&3 \\
   \hline   & $\Bbb Z_2[X,Y]/(X^3,XY,Y^2) $ & 16&  Fig. 5  & 4 & 1&4\\
    \hline  & $\Bbb Z_8[X]/(2X,X^2)$ & 16 &   Fig. 5  & 4 & 1&4 \\
    \hline  & $\Bbb Z_4[X]/(X^3,2X^2,2X)$ & 16 &   Fig. 5  & 4 & 1&4 \\
 \hline   & $\Bbb Z_4[X]/(X^2+2X) $& 16 &  Fig. 6  & 3& 1&3 \\
  \hline  & $\Bbb Z_8[X]/(2X,X^2+4)$ & 16 &   Fig. 6  & 3& 1&3 \\
   \hline  & $\Bbb  Z_2[X,Y]/(X^2,Y^2-XY)$ & 16 &  Fig. 6 & 3& 1&3 \\
   \hline   & $\Bbb Z_4[X,Y]/(X^2,Y^2-XY,XY-2,2X,2Y) $ & 16&  Fig. 6 & 3& 1&3\\
 \hline   & $\Bbb  Z_4[X,Y]/(X^2,Y^2,XY-2,2X,2Y) $& 16 &  Fig. 7 & 3& 1&3 \\
  \hline  & $\Bbb  Z_2[X,Y]/(X^2,Y^2) $ & 16 &  Fig. 7 & 3& 1&3\\
   \hline  & $\Bbb  Z_4[X]/(X^2)$ & 16 & Fig. 7 & 3& 1&3 \\
    \hline   & $\Bbb Z_4[X]/(X^3-X^2-2,2X^2,2X) $ & 16& Fig. 8 & 4 & 1 &3\\
 \hline   & $\Bbb  Z_2[X,Y,Z]/(X,Y,Z)^2 $& 16 & $K_7$ & 1& 1&7 \\
  \hline  & $\Bbb  Z_4[X,Y]/(X^2,Y^2,XY,2X,2Y) $ & 16 &  $K_7$ &1& 1&7 \\
   \hline  & $\Bbb F_8[X]/(X^2)$ & 64 &  $K_7$ & 1& 1&7\\
   \hline  & $\Bbb  Z_4[X]/(X^3+X+1)$ & 64 & $K_7$ & 1& 1&7 \\
    \hline
\end{tabular}
\]




\begin{figure}[h]
\begin{minipage}{4.8cm}
\includegraphics[width=\textwidth]{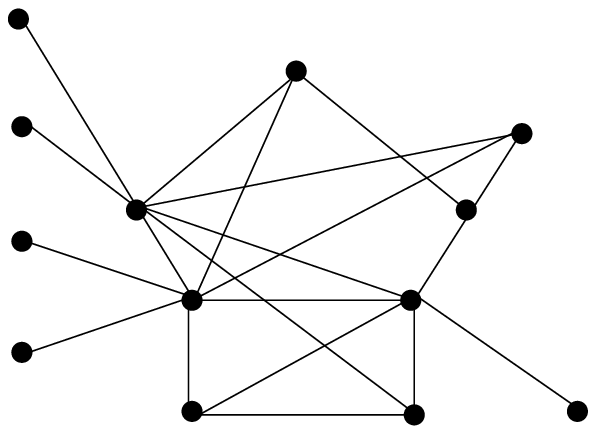}
\caption{\label{figure19} Graph for $\mathbb{Z}_2\times \mathbb{Z}_2 \times \mathbb{Z}_4$ and $\mathbb{Z}_2\times \mathbb{Z}_2 \times \mathbb{Z}_2[X]/(X^2) $ }
\end{minipage}
\begin{minipage}{4.8cm}
\includegraphics[width=4.5cm,height=2.92cm]{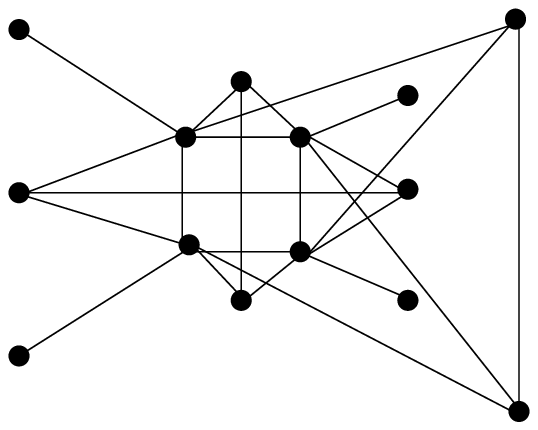}
\caption{\label{figure20} Graph for $\mathbb{Z}_2\times \mathbb{Z}_2 \times \mathbb{Z}_2\times \mathbb{Z}_2 $}
\end{minipage}
\begin{minipage}{4.8cm}
\includegraphics[width=4.5cm,height=3.5cm]{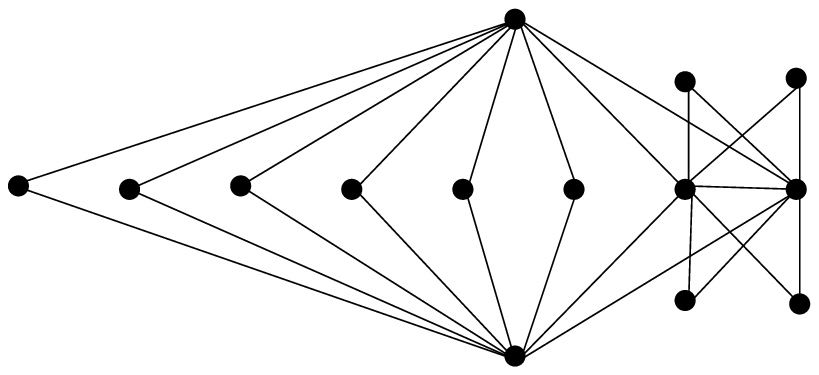}
\caption{\label{figure21} Graph for Graph for $\mathbb{Z}_3\times \mathbb{Z}_9$ and $\mathbb{Z}_3\times \mathbb{Z}_3[X]/(X^2) $.}
\end{minipage}
\end{figure}

\[
\begin{tabular}{|c|c|c|c|c|c|c|}
  \hline
  Vertices & R &$ |R|$ & Graph &   $\alpha(\Gamma(R))$ & $\gamma(\Gamma(R))$&$\omega(\Gamma(R)) $\\
 \hline 8 & $\Bbb Z_2\times\Bbb  F_8 $& 16 & $K_{1,7}$ & 7 & 1&2 \\
 \hline   & $\Bbb Z_3\times\Bbb Z_7 $ & 21& $K_{2,6}$ & 6 & 2 &2\\
 \hline   & $\Bbb Z_5\times\Bbb Z_5 $&25 & $K_{4,4}$ & 4& 2&2 \\
  \hline  & $\Bbb Z_{27}$ & 27 &  Fig. 9 & 6 & 1&3 \\
   \hline  & $\Bbb  Z_9[X]/(3X,X^2-3)$ & 27 & Fig. 9 &  6 & 1&3 \\
    \hline  & $\Bbb  Z_9[X]/(3X,X^2-6)$ & 27 & Fig. 9 &  6 & 1&3 \\
 \hline  & $\Bbb Z_3[X]/(X^3) $& 27 & Fig. 9 &  6 & 1&3 \\
 \hline   & $\Bbb Z_3[X,Y]/(X,Y)^2 $ & 27& $K_8$ & 1 & 1 &8\\
 \hline   & $\Bbb Z_9[X]/(3,X)^2 $& 27 & $K_8$ & 1& 1&8 \\
  \hline  & $\Bbb F_9[X]/(X^2)$ & 81 &  $K_8$ & 1 & 1&8 \\
   \hline  & $\Bbb  Z_9[X]/(X^2+1)$ & 81 &$K_8$ & 1 & 1&8 \\
    \hline
\end{tabular}
\]

\[
\begin{tabular}{|c|c|c|c|c|c|c|}
  \hline
  Vertices & R &$ |R|$ & Graph &   $\alpha(\Gamma(R))$ & $\gamma(\Gamma(R))$&$\omega(\Gamma(R)) $\\
 \hline 9 & $\Bbb Z_2\times   F_9 $& 18 & $K_{1,8}$ & 8 & 1&2 \\
 \hline   & $\Bbb Z_3\times  F_8 $ & 24& $K_{2,7}$ & 7 & 2 &2\\
 \hline   & $ F_4\times\Bbb Z_7 $& 28 & $K_{3,6}$ & 6& 2&2 \\
  \hline  & $\Bbb Z_2\times\Bbb Z_2\times\Bbb Z_3 $ & 12 &  Fig. 10 & 6 & 3&3 \\
  \hline   & $\Bbb Z_4\times F_4 $& 16 & Fig. 11 & 6& 2&2 \\
   \hline  & $\Bbb  Z_2[X]/(X^2)\times  F_4$ & 16 & Fig. 11 & 6 & 2&2 \\
    \hline
\end{tabular}
\]

\[
\begin{tabular}{|c|c|c|c|c|c|c|}
  \hline
  Vertices & R &$ |R|$ & Graph &   $\alpha(\Gamma(R))$ & $\gamma(\Gamma(R))$&$\omega(\Gamma(R)) $\\
 \hline 10 & $\Bbb Z_3\times   F_9 $& 27 & $K_{2,8}$ & 8 & 2&2 \\
 \hline   & $\Bbb F_4\times  F_8 $ & 32& $K_{3,7}$ & 7 & 2 &2\\
 \hline   & $\Bbb Z_5\times\Bbb Z_7 $& 35 & $K_{4,6}$ & 6& 2&2 \\
  \hline  & $\Bbb Z_{121} $ & 121 & $K_{10}$ & 1 & 1&10 \\
  \hline   & $\Bbb Z_{11}[X]/(X^2) $& 121 & $K_{10}$ & 1& 1&10 \\
    \hline
\end{tabular}
\]

\[
\begin{tabular}{|c|c|c|c|c|c|c|}
  \hline
  Vertices & R &$ |R|$ & Graph &   $\alpha(\Gamma(R))$ & $\gamma(\Gamma(R))$&$\omega(\Gamma(R)) $\\
 \hline 11 & $\Bbb Z_2\times\Bbb  Z_{11} $& 22 & $K_{1,10}$ & 10 & 1&2 \\
 \hline   & $ F_4\times\Bbb F_9 $ & 36& $K_{3,8}$ & 8 & 2 &2\\
 \hline   & $\Bbb Z_5\times  F_8 $& 40 & $K_{4,7}$ & 7& 2&2 \\
  \hline  & $\Bbb Z_2\times\Bbb Z_9 $ & 18 &  Fig. 12 & 8 & 3&3 \\
  \hline   & $\Bbb Z_2\times\Bbb Z_3[X]/(X^2) $& 18 & Fig. 12 & 8& 3&3\\
   \hline  & $\Bbb  Z_5\times \Bbb Z_4$ & 20 & Fig. 13 & 8 & 3&2 \\
  \hline  & $\Bbb Z_5\times\Bbb Z_2[X]/(X^2) $ & 20 &  Fig. 13 & 8 & 3&2 \\
  \hline   & $\Bbb Z_2\times\Bbb Z_8 $& 16 & Fig. 14 & 8& 2&3 \\
   \hline  & $\Bbb  Z_2\times\Bbb Z_2[X]/(X^3)$ & 16 & Fig. 14 &  8& 2&3  \\
  \hline  & $\Bbb Z_2\times\Bbb Z_4[X]/(2X,X^2-2) $ & 16 &  Fig. 14 &  8& 2&3  \\
  \hline   & $\Bbb Z_2\times\Bbb Z_2[X,Y]/(X,Y)^2 $& 16 & Fig. 15 & 7& 2&4 \\
   \hline  & $\Bbb  Z_2\times\Bbb Z_4[X]/(2,X)^2$ & 16 & Fig. 15 & 7 & 2&4 \\
  \hline  & $\Bbb Z_4\times\Bbb Z_4 $ & 16 &  Fig. 16 & 6 & 2&3 \\
  \hline   & $\Bbb Z_4\times\Bbb Z_2[X]/(X^2) $& 16 & Fig. 16 & 6 & 2&3\\
   \hline  & $\Bbb Z_2[X]/(X^2)\times \Bbb Z_2[X]/(X^2)$ & 16 & Fig. 16 & 6 & 2&3 \\
    \hline
\end{tabular}
\]

\[
\begin{tabular}{|c|c|c|c|c|c|c|}
  \hline
  Vertices & R &$ |R|$ & Graph &   $\alpha(\Gamma(R))$ & $\gamma(\Gamma(R))$&$\omega(\Gamma(R)) $\\
 \hline 12 & $\Bbb Z_3\times\Bbb  Z_{11} $& 33 & $K_{2,10}$ & 10 & 2&2 \\
 \hline   & $\Bbb Z_5\times\Bbb Z_9 $ & 45& $K_{4,8}$ & 8 & 2 &2\\
 \hline   & $\Bbb Z_7\times\Bbb Z_7 $& 49 & $K_{6,6}$ & 6& 2&2 \\
 \hline   & $\Bbb Z_2\times\Bbb Z_2\times\Bbb Z_4 $& 16 & Fig. 17 & 6& 2&2 \\
  \hline  & $\Bbb Z_{169} $ & 169 & $K_{12}$ & 1 & 1&12 \\
  \hline   & $\Bbb Z_{13}[X]/(X^2) $& 169 & $K_{12}$ & 1& 1&12 \\
    \hline
\end{tabular}
\]

\[
\begin{tabular}{|c|c|c|c|c|c|c|}
  \hline
  Vertices & R &$ |R|$ & Graph &   $\alpha(\Gamma(R))$ & $\gamma(\Gamma(R))$&$\omega(\Gamma(R)) $\\
 \hline 13 & $\Bbb Z_2\times\Bbb  Z_{13} $& 26 & $K_{1,12}$ & 12 & 1&2 \\
 \hline   & $ F_4\times\Bbb Z_{11} $ & 44& $K_{3,10}$ & 10 & 2 &2\\
 \hline   & $\Bbb Z_7\times  F_8 $& 56 & $K_{6,7}$ & 7& 2&2 \\
  \hline  & $\Bbb  Z_2\times\Bbb Z_3\times\Bbb Z_3 $  & 18 & Fig. 18& 8 & 3&3 \\
  \hline  & $\Bbb  Z_2\times\Bbb Z_2\times\Bbb Z_4 $  & 16 & Fig. 19 & 8 & 3&3 \\
  \hline   & $\Bbb  Z_2\times\Bbb Z_2\times\Bbb Z_2[X]/(X^2) $& 16 & Fig. 19 & 8& 3&3 \\
    \hline
\end{tabular}
\]

\[
\begin{tabular}{|c|c|c|c|c|c|c|}
  \hline
  Vertices & R &$ |R|$ & Graph &   $\alpha(\Gamma(R))$ & $\gamma(\Gamma(R))$&$\omega(\Gamma(R)) $\\
 \hline 14 & $\Bbb Z_3\times\Bbb  Z_{13} $& 39 & $K_{2,12}$ & 12 & 2&2 \\
 \hline   & $\Bbb Z_5\times\Bbb Z_{11} $ & 55& $K_{4,10}$ & 10 & 2 &2\\
 \hline   & $\Bbb Z_7\times  F_9 $& 63 & $K_{6,8}$ & 8& 2&2 \\
  \hline  & $\Bbb  F_8\times  F_8 $  & 64 &$K_{7,7}$& 7 & 2&2 \\
  \hline  & $\Bbb  Z_2\times\Bbb Z_2\times Z_2\times\Bbb Z_2 $  & 16 & Fig. 20 & 7 & 4&3 \\
  \hline   & $\Bbb  Z_3\times\Bbb Z_9 $& 27 & Fig. 21 & 10& 2&3 \\
    \hline   & $\Bbb Z_3\times\Bbb Z_3[X]/(X^2) $& 27 & Fig. 21 & 10& 2&3 \\
    \hline
\end{tabular}
\]



\end{document}